\documentclass[12pt]{amsart}
\usepackage{amssymb,latexsym,comment,url}

\newcommand{\LL}{{\mathcal L}}
\newcommand{\EE}{{\mathcal E}}
\newcommand{\Cdes}{\widetilde{\mathcal C}}
\newcommand{\Char}{\operatorname{char}}
\newcommand{\Diag}{\operatorname{Diag}}
\newcommand{\Spec}{\operatorname{Spec}}
\newcommand{\la}{\lambda}
\newcommand{\GL}{\operatorname{GL}}
\newcommand{\divv}{\operatorname{div}}
\newcommand{\K}{{\mathbb K}}
\newcommand{\Kbar}{\overline{K}}
\newcommand{\isom}{ \cong }
\newcommand{\nr}{{\text{\rm nr}}}
\newcommand{\OK}{{\mathcal{O}_K}}
\newcommand{\pf}{\operatorname{pf}}
\newcommand{\Pf}{\operatorname{Pf}}
\newcommand{\CC}{\operatorname{\mathcal C}}
\newcommand{\PP}{{\mathbb P}}
\newcommand{\Q}{{\mathbb Q}}
\newcommand{\OO}{{\mathcal O}}
\newcommand{\Sing}{\operatorname{Sing}}
\newcommand{\rank}{\operatorname{rank}}
\newcommand{\ra}{{\longrightarrow}}
\newcommand{\Z}{{\mathbb Z}}
\newcommand{\adj}{\operatorname{adj}}

\newenvironment{Proof}{\par\noindent{\sc Proof:}}%
                      {\hspace*{\fill}\nobreak$\Box$\par\medskip}
\newenvironment{ProofOf}[1]{\par\noindent{\sc Proof of #1:}}%
                       {\hspace*{\fill}\nobreak$\Box$\par\medskip}

\newtheorem{Proposition}{Proposition}[section]
\newtheorem{Theorem}[Proposition]{Theorem}
\newtheorem{Lemma}[Proposition]{Lemma}
\newtheorem{Corollary}[Proposition]{Corollary}

\theoremstyle{definition}
\newtheorem{Definition}[Proposition]{Definition}
\newtheorem{Remark}[Proposition]{Remark}

\renewcommand{\baselinestretch}{1.1}

\begin{document}

\title[Genus one curves of degree 5]%
{On genus one curves of degree 5 with square-free discriminant}

\author{Tom~Fisher}
\address{University of Cambridge,
         DPMMS, Centre for Mathematical Sciences,
         Wilberforce Road, Cambridge CB3 0WB, UK}
\email{T.A.Fisher@dpmms.cam.ac.uk}

\author{Mohammad~Sadek}
\address{American University in Cairo, Mathematics and Actuarial Science Department, AUC Avenue, New Cairo, Egypt}
\email{mmsadek@aucegypt.edu}

\date{26th June 2014}

\begin{abstract}
  We study genus one curves of degree $5$ defined by Pfaffians.  We
  give new formulae for the invariants, and prove the equivalence of
  two different definitions of minimality. As an application we show
  that transformations between models with square-free discriminant
  are necessarily integral. This result is used by Bhargava and
  Shankar in their work on the average ranks of elliptic curves.
\end{abstract}

\maketitle


\section{Introduction}

Let $E$ be an elliptic curve over a number field $\K$. An {\em
  $n$-covering} of $E$ is a pair $(C,\pi)$ where $C$ is a smooth curve
of genus one and $\pi:C\to E$ is a morphism such that $\pi=[n]\circ
\psi$ for some isomorphism $\psi :C\to E$ defined over
$\overline{\K}$. If $C$ is everywhere locally soluble then by
\cite[Theorem 1.3]{CaIV} there exists a $\K$-rational divisor $D$ on
$C$ such that $D$ is linearly equivalent to $\psi^*(n.0_E)$. The
linear system $|D|$ defines a morphism $C\to \PP^{n-1}$. If $n \ge 3$
then this morphism is an embedding, and the image is called a {\em
  genus one normal curve} of degree $n$.  The word ``normal'' refers
to the fact the curve is projectively normal, i.e. the homogeneous
co-ordinate ring is integrally closed.  This should not be confused
with the fact $C$ is normal, which is automatic since $C$ is smooth.

When $n=2,3,4$ the curve $C$ is represented by a binary quartic,
ternary cubic, or pair of quadrics in 4 variables.  In this paper we
take $n=5$, in which case $C$ is represented by data of the following
form.

A {\em Pfaffian model} $\Phi$ over a ring $R$ is a $5 \times 5$
alternating matrix of linear forms in $R[x_1, \ldots,x_5]$.  We write
$X_5(R)$ for the space of all Pfaffian models over $R$.  Two models
$\Phi$ and $\Phi'$ are $R$-equivalent if $\Phi'=[A,B]\Phi$ for some
$A,B\in \GL_5(R)$. The action of $A$ is given by $\displaystyle
\Phi\mapsto A\Phi A^T$, and the action of $B$ is given by \[
(\Phi_{ij}(x_1,\dots,x_5))\mapsto (\Phi_{ij}(x_1',\ldots,x_5')) \]
where $x_j'=\sum_{i=1}^5B_{ij}x_i$. We define $\det [A,B] =(\det
A)^2\det B$.  The models $\Phi$ and $\Phi'$ are {\em properly
  $R$-equivalent} if $\det [A,B] = 1$.  The invariants $c_4, c_6,
\Delta \in \Z[X_5]$ are certain integer coefficient polynomials in the
$50$ coefficients of a Pfaffian model.  We give formulae for these in
Section~\ref{sec:inv}.

We work over a discrete valuation field $K$ with valuation ring $\OK$,
normalised valuation $v : K^\times \to \Z$, uniformiser $\pi$, and
residue field $k = \OK/\pi \OK$.

Our main result is the following.  It answers a question of Bhargava,
and is used in the work of Bhargava and Shankar~\cite[Proposition
11]{BS5} on the average size of the $5$-Selmer group of an elliptic
curve.

\begin{Theorem}
\label{thm1}
Let $\Phi, \Phi' \in X_5(\OK)$ be Pfaffian models with
$v(\Delta(\Phi)) \le 1$ and $v(\Delta(\Phi')) \le 1$.  If $\Phi' =
[A,B] \Phi$ for some $A,B \in \GL_5(K)$ then $A,B \in K^\times
\GL_5(\OK)$.  In particular
\begin{enumerate}
\item If $\Phi$ and $\Phi'$ are $K$-equivalent then they are
  $\OK$-equivalent.
\item The stabiliser of $\Phi$ in $\GL_5(K) \times \GL_5(K)$ is
  contained in the subgroup generated by $\GL_5(\OK) \times
  \GL_5(\OK)$ and $[\pi^{-1} I_5, \pi^2 I_5]$.
\end{enumerate}
\end{Theorem}

To indicate how Theorem~\ref{thm1} is useful, we give the following
global application. We take $\K=\Q$, but note that the result
generalises immediately to any number field with class number $1$.  We
say that a Pfaffian model $\Phi$ has the same invariants as an
elliptic curve $E$ if the invariants $c_4(\Phi)$, $c_6(\Phi)$,
$\Delta(\Phi)$ are the same as the invariants $c_4, c_6, \Delta$ of a
minimal Weierstrass equation for $E$.

\begin{Theorem}
\label{thm:glob}
Let $E/\Q$ be an elliptic curve with square-free minimal
discriminant. Then the $5$-Selmer group $S^{(5)}(E/\Q)$ is in
bijection with the set of Pfaffian models over $\Z$ with the same
invariants as $E$, up to proper $\Z$-equivalence.
\end{Theorem}

In Sections~\ref{sec:min} and~\ref{sec:geom-min} we introduce two
different definitions of minimality, and show that if they agree then
Theorem~\ref{thm1} is a natural consequence.  The agreement of the two
definitions is proved in Sections~\ref{sec:flat}, \ref{sec:normal} and
\ref{sec:invdiff}. This extends \cite[Theorem~4.1]{minimal} from genus
one curves of degrees $2$, $3$ and $4$, to degree $5$.  In
Section~\ref{sec:alt} we give a short alternative proof of
Theorem~\ref{thm1}, that is motivated by the ideas in the rest of this
paper, but avoids nearly all the scheme-theoretic machinery.


\section{Pfaffians and invariants}
\label{sec:inv}

In this section we briefly describe how the equations for a genus one
normal curve of degree~$5$ can be written in terms of Pfaffians.  We
then give some new formulae for the invariants of a Pfaffian model,
that are simpler than the evaluation algorithms in \cite[Section
8]{g1inv}.

The Pfaffian of an alternating matrix is an integer coefficient
polynomial in the entries of the matrix, whose square is the
determinant. We only need to consider Pfaffians of $4 \times 4$
matrices, in which case
\begin{equation*}
 \pf \begin{pmatrix} 0 & a_{12} & a_{13} & a_{14} \\
& 0 & a_{23} & a_{24} \\
& & 0 & a_{34} \\
& & & 0
\end{pmatrix} = a_{12} a_{34} - a_{13} a_{24} + a_{14} a_{23}.
\end{equation*}
If $\Phi$ is an $5\times 5$ alternating matrix then the row vector of
submaximal Pfaffians of $\Phi$ is $\Pf(\Phi)=(p_1,\ldots,p_5)$ where
$p_i=(-1)^{i}\pf\left(\Phi^{\{i\}}\right)$ and $\Phi^{\{i\}}$ is the
matrix obtained by deleting the $i$th row and column of $\Phi$.  It
can be shown, for example by direct calculation, that
$\Pf(\Phi)\Phi=0$, $\adj(\Phi)=\Pf(\Phi)^T\Pf(\Phi)$ and $\Pf(A\Phi
A^T)=\Pf(\Phi)\adj(A)$ for all $5 \times 5$ matrices $A$.

In this section we work over any field $K$.  Let $C \subset \PP^4_K$
be a genus one normal curve, i.e. a smooth curve of genus one embedded
by a complete linear system of degree $5$.  Let $R = K[x_1, \ldots,
x_5] = \oplus_{d \ge 0} R_d$ be the polynomial ring with its usual
grading by degree. Let $R(d)$ be the graded free $R$-module of rank
$1$ with $R(d)_e = R_{d+e}$.  By the Buchsbaum-Eisenbud structure
theorem~\cite{BE1}, \cite{BE2}, or the treatment specific to this case
in \cite{explicit5}, the coordinate ring of $C$ has minimal free
resolution
\begin{equation}
\label{mfr0}
 0 \ra R(-5) \stackrel{\Pf(\Phi)^T}{\ra} R(-3)^5 \stackrel{\Phi}{\ra}
R(-2)^5 \stackrel{\Pf(\Phi)}{\ra} R
\end{equation}
for some $\Phi \in X_5(K)$.  In particular the homogeneous ideal of
$C$ is generated by the $4 \times 4$ Pfaffians of $\Phi$.  More
generally, for any $\Phi \in X_5(K)$, we let $C_\Phi \subset \PP^4_K$
be the subscheme defined by its $4 \times 4$ Pfaffians. We say that
$\Phi$ is {\em non-singular} if $C_\Phi$ is a smooth curve of genus
one.  We write $K[X_5]$ for the polynomial ring in the $50$
coefficients of a Pfaffian model. A polynomial $F \in K[X_5]$ is an
{\em invariant of weight $k$} if $F \circ g=(\det g)^k F$ for all $g
\in \GL_5 \times \GL_5$.

\begin{Theorem}\label{thm-inv}
  There are invariants $c_4,c_6,\Delta\in\Z[X_5]$ of degrees
  $20,30,60$ and weights $4,6,12$, satisfying
  $c_4^3-c_6^2=1728\Delta$, with the following properties.
\begin{enumerate}
\item If $\Char(K) \not=2,3$ then the ring of invariants in $K[X_5]$
  is generated by (the images of) $c_4$ and $c_6$.
\item A model $\Phi\in X_5(K)$ is non-singular if and only if
  $\Delta(\Phi)\ne 0$.
\item There exist $a_1,a_2,a_3,a_4,a_6,b_2,b_4,b_6 \in \Z[X_5]$
  satisfying
    \begin{equation*}
      \begin{aligned}
        b_2 & = a_1^2+4 a_2, & b_4 & = a_1 a_3 + 2 a_4,
        \quad \quad b_6  =  a_3^2 + 4 a_6, \\
        c_4 & = b_2^2-24 b_4, \quad & c_6 & = -b_2^3 + 36 b_2 b_4 -
        216 b_6,
      \end{aligned}
    \end{equation*}
    such that if $\Phi \in X_5(K)$ is non-singular then $\CC_\Phi$ has
    Jacobian
    \begin{equation}
      \label{eqn:jac}
      y^2 + a_1(\Phi) xy + a_3(\Phi) y
      = x^3 + a_2(\Phi) x^2 + a_4(\Phi) x + a_6(\Phi).
    \end{equation}
\end{enumerate}
\end{Theorem}
\begin{Proof} This is \cite[Theorem 4.4]{g1inv} together with
  \cite[Theorem 1.1]{minred5}.
\end{Proof}

It is shown in \cite[Section 5.4]{g1inv} that if $\Char K \not=2$ and
$\Phi \in X_5(K)$ is non-singular then there is an invariant
differential $\omega_\Phi$ on $C_\Phi$ given by
\begin{equation}
\label{invdiff}
\omega_{\Phi}=\frac{x_i^2d(x_j/x_i)}{Q(x_1,\ldots,x_5)},
\textrm{ where }
Q =\frac{\partial P}{\partial x_k}
\frac{\partial \Phi}{\partial x_l}\frac{\partial P^T}{\partial x_m},
\end{equation}
$P = \Pf(\Phi)$, and $(i,j,k,l,m)$ is any even permutation of
$(1,2,3,4,5)$. In the definition of $Q$, it is understood that by the
partial derivative of a matrix we mean the matrix of partial
derivatives.  As we show in Remark~\ref{char2}, the restriction $\Char
K \not=2$ is not needed.

In \cite[Section 7]{invenqI} an alternative description of the
invariant differential is given in terms of a certain covariant.  We
now give an explicit construction of this covariant, based in part on
ideas in \cite[Section 4]{B5}.  For $(i,j,k,l,m)$ an even permutation
of $(1,2,3,4,5)$ we define
\[ \Omega_{ij} = \frac{\partial P}{\partial x_k} \frac{\partial
  \Phi}{\partial x_l} \frac{\partial P^T}{\partial x_m} +
\frac{\partial P}{\partial x_m} \frac{\partial \Phi}{\partial x_k}
\frac{\partial P^T}{\partial x_l} + \frac{\partial P}{\partial x_l}
\frac{\partial \Phi}{\partial x_m} \frac{\partial P^T}{\partial x_k}.
\]
Now $\Omega = (\Omega_{ij})$ is an alternating matrix of quadratic
forms. We define an action of $\GL_5 \times \GL_5$ on the space of
such matrices via
\[ [A,B] : \Omega \mapsto B^{-T} (\Omega_{ij}(x_1',\ldots,x_5'))
B^{-1} \] where $x_j'=\sum_{i=1}^5B_{ij}x_i$. In particular the first
copy of $\GL_5$ acts trivially. Recall that for $g = [A,B]$ we defined
$\det g = (\det A)^2 \det B$.

\begin{Lemma}
\label{omega-cov}
The map $\Phi \mapsto \Omega$ is a covariant of weight $1$, in the
sense that
\[ g \Phi \mapsto (\det g) g \Omega \] for all $g \in \GL_5 \times
\GL_5$.
\end{Lemma}
\begin{Proof}
  If we replace $\Phi$ by $A \Phi A^T$ then $P$ is replaced by $P \adj
  A$ and $\Omega$ is multiplied by $(\det A)^2$. So it suffices to
  consider $g = [I_5,B]$ for $B$ running over a set of generators for
  $\GL_5$. Since the cases where $B$ is a diagonal matrix or a
  permutation matrix are easy, this reduces us to considering $B = I_5
  + \lambda E_{12}$, where $\lambda \in K$ and $E_{12}$ is the
  elementary matrix with a $1$ in position $(1,2)$ and all other
  entries $0$. This corresponds to the substitution $x_2 \leftarrow
  x_2 + \la x_1$. In the definition of $\Omega_{ij}$ we replace
  $\frac{\partial P}{\partial x_1}$ by $\frac{\partial P}{\partial
    x_1} + \la \frac{\partial P}{\partial x_2}$ and $\frac{\partial
    \Phi}{\partial x_1}$ by $\frac{\partial \Phi}{\partial x_1} + \la
  \frac{\partial \Phi}{\partial x_2}$.  This has the effect of
  replacing $\Omega_{r2}$ by $\Omega_{r2} - \la \Omega_{r1}$ and
  $\Omega_{2r}$ by $\Omega_{2r} - \la \Omega_{1r}$ for $r=3,4,5$.  A
  calculation, using the fact $\Phi$ is alternating, shows that the
  other entries of $\Omega$ do not change. Thus $\Omega$ changes to $g
  \Omega$ as required.
\end{Proof}

We put
\[ M_{ij} = \sum_{r,s=1}^5 \frac{\partial \Omega_{ir}}{\partial x_s}
\frac{\partial \Omega_{js}}{\partial x_r} \quad \text{ and } \quad
N_{ijk} = \sum_{r=1}^5 \frac{\partial M_{ij}}{\partial x_r}
\Omega_{rk}. \]

\begin{Theorem}
The invariants $c_4,c_6 \in \Z[X_5]$ are given by
\[ c_4(\Phi) =
\frac{1}{13440} 
\sum_{i,j,r,s=1}^5 \frac{ \partial^2 M_{ij}}{\partial x_r \partial
  x_s} \frac{ \partial^2 M_{rs}}{\partial x_i \partial x_j} \] and
\[ c_6(\Phi) =
\frac{-1}{1036800} 
\sum_{i,j,k,r,s,t=1}^5 \frac{ \partial^3 N_{ijk}}{\partial
  x_r \partial x_s
\partial x_t }
 \frac{\partial^3 N_{rst}}{\partial x_i \partial x_j \partial x_k}. \]
\end{Theorem}
\begin{Proof}
  It may be checked using Lemma~\ref{omega-cov} that these polynomials
  are invariants of degrees $20$ and $30$.  By Theorem~\ref{thm-inv}
  it only remains to show they are scaled as specified in
  \cite{g1inv}. We can do this by computing a single numerical
  example.
\end{Proof}

We may compute the discriminant $\Delta$ either as $(c_4^3 -
c_6^2)/1728$, or directly using the method at the end of \cite[Section
8]{g1inv}.

\section{Minimal Pfaffian models}
\label{sec:min}

In this section we make some remarks about minimal Pfaffian models,
and more specifically those with square-free discriminant.  We also
explain how Theorem~\ref{thm:glob} follows from Theorem~\ref{thm1}.

From now on $K$ will be a discrete valuation field, with ring of
integers $\OK$, and normalised valuation $v:K^{\times}\to\Z$.  We fix
a uniformiser $\pi$ and write $k$ for the residue field.  Let
$S=\Spec\OK$. For the proof of Theorem~\ref{thm1} we are free to
replace $K$ by any unramified extension. We may therefore assume when
convenient that $K$ is complete, and $k$ is algebraically closed.

A Pfaffian model $\Phi \in X_5(K)$ is {\em integral} if $\Phi \in
X_5(\OK)$, i.e. it has coefficients in $\OK$.  It follows from
Theorem~\ref{thm-inv} that if $\Phi$ is non-singular and integral then
$v(\Delta(\Phi))=v(\Delta_E)+12 \ell(\Phi)$ where $\Delta_E$ is the
minimal discriminant of the Jacobian $E$, and $\ell(\Phi) \ge 0$ is an
integer called the {\em level}.  We say that $\Phi$ is {\em minimal}
if $v(\Delta(\Phi))$ is minimal among all integral models
$K$-equivalent to $\Phi$. If $\Phi'=g \Phi$ for $g=[A,B]$ with $A,B\in
\GL_5(K)$ than $\ell(\Phi')=\ell(\Phi)+v(\det g)$.

\begin{Theorem}
\label{thm2}
(Minimisation theorem)
Let $\Phi\in X_5(K)$ be non-singular.
If $C_{\Phi}(K)\ne \emptyset$ then $\Phi$ is $K$-equivalent
to an integral model of level $0$.
\end{Theorem}
\begin{Proof}
This is \cite[Theorem 2.1(i)]{minred5}.
\end{Proof}

The proof of Theorem~\ref{thm2} is rather short.
In \cite{minred5} the first author also investigated to what
extent the hypothesis $C_{\Phi}(K)\ne \emptyset$ can be weakened,
and gave an algorithm for minimising.

\begin{Lemma}
\label{lem:sqrfree}
Let $\Phi \in X_5(\OK)$ with $v(\Delta(\Phi)) \le 1$.
\begin{enumerate}
\item The Jacobian $E$ of $C_\Phi$ has Kodaira symbol $I_0$ or $I_1$.
\item If $K$ is a $p$-adic field then $C_{\Phi}(K)\ne \emptyset$.
\end{enumerate}
\end{Lemma}
\begin{Proof} (i) By Theorem~\ref{thm-inv} we have $v(\Delta_E) \le
  1$.  It follows by Tate's algorithm that the Kodaira symbol is
  either
  $I_0$ or $I_1$. \\
  (ii) Since $v(\Delta(\Phi)) < 12$ we have $\ell(\Phi)=0$.  Then by
  \cite[Theorem 7.1]{minred5} we have $C_{\Phi}(K^\nr)\ne \emptyset$
  where $K^\nr$ is the maximal unramified extension. By (i) we know
  that $E/K$ has Tamagawa number $1$. Therefore, as explained in
  \cite[Lemma 2.1]{FS}, solubility over $K^\nr$ is equivalent to
  solubility over $K$.
\end{Proof}

\begin{Remark}
\label{rem1}
To prove Theorem~\ref{thm1} it suffices to show that $B \in K^\times
\GL_5(\OK)$. The reason for this is as follows.  By
Lemma~\ref{lem:i-ii} we know that if $\Phi$ is minimal then its $4
\times 4$ Pfaffians are linearly independent mod $\pi$.  So if $\Phi$
and $\Phi'$ are both minimal and $\Phi' = [A,\la I_5] \Phi$, then from
$\Pf(\Phi') = \la^2 \Pf(\Phi) \adj(A)$ we deduce that $A \in K^\times
\GL_5(\OK)$. The final statements (i) and (ii) of Theorem~\ref{thm1}
are immediate, since $v(\det[A,B]) = 0$ and the transformations $[\la
I_5,\la^{-2} I_5]$ for $\la \in K^\times$ act trivially on the space
of Pfaffian models.
\end{Remark}

We now explain how Theorem~\ref{thm:glob} follows from
Theorem~\ref{thm1}.

\begin{Theorem}
\label{thm:sel5}
Let $E/\Q$ be an elliptic curve. The $5$-Selmer group $S^{(5)}(E/\Q)$
is in bijection with the proper $\Q$-equivalence classes of Pfaffian
models $\Phi \in X_5(\Q)$ with the same invariants as $E$ and
$C_\Phi(\Q_p) \not= \emptyset$ for all primes $p$.
\end{Theorem}
\begin{Proof}
This is a special case of \cite[Theorem 6.1]{invenqI}.
\end{Proof}

\begin{ProofOf}{Theorem~\ref{thm:glob}}
  By Theorem~\ref{thm2} and strong approximation, each of the classes
  in Theorem~\ref{thm:sel5} contains a model with coefficients in
  $\Z$. Since $\Delta_E$ is square-free, Theorem~\ref{thm1} shows that
  the map from proper $\Z$-equivalence classes to proper
  $\Q$-equivalence classes is injective.  Moreover the condition
  $C_\Phi(\Q_p) \not= \emptyset$ is automatically satisfied by
  Lemma~\ref{lem:sqrfree}.
\end{ProofOf}

Let $\Phi \in X_5(\OK)$ have reduction $\phi \in X_5(k)$.  We write
$\CC_\Phi \subset \PP^4_S$ for the $S$-scheme defined by the $4 \times
4$ Pfaffians. It has generic fibre $C_\Phi$ and special fibre
$C_\phi$.

Suppose the entries of $\phi$ span $\langle x_1, \ldots, x_5 \rangle$.
If $P$ is $k$-point on $C_\phi$ then by an $\OK$-equivalence we may
assume $P = (1:0:\ldots:0)$.  We may further assume $\phi_{12} = x_1$
and all other $\phi_{ij}$ (for $i < j$) are linear forms in $x_2,
\ldots ,x_5$.  The tangent space to $C_\phi$ at $P$ is $\{\phi_{34} =
\phi_{35} = \phi_{45} = 0\} \subset \PP_k^4$.

\begin{Lemma} 
\label{lem:regdef}
Let $P \in C_\phi$ as above. The following are equivalent.
\begin{enumerate}
\item The tangent space to $\CC_\Phi$ at $P$ has dimension at most 2.
\item Every linear combination $r \Phi_{34} + s \Phi_{35} + t \Phi_{45}$ 
(where $r,s,t \in \OK$, not all in $\pi \OK$) that vanishes mod $\pi$
has coefficient of $x_1$ not divisible by $\pi^2$.
\end{enumerate}
\end{Lemma}
\begin{Proof} By (i) we mean $\dim ({\mathfrak m}_P/{\mathfrak m}_P^2)
  \le 2$ where ${\mathfrak m}_P$ is the maximal ideal of the local
  ring at $P$. The lemma is proved by a straightforward calculation.
\end{Proof}

The following lemma will be used both to show that $\CC_\Phi$ is
regular, and in the elementary proof of Theorem~\ref{thm1} in
Section~\ref{sec:alt}.

\begin{Lemma}
\label{lem:reg}
If $\Phi \in X_5(\OK)$ with $v(\Delta(\Phi)) \le 1$ then every
$k$-point $P$ on $C_\phi$ satisfies the conditions in
Lemma~\ref{lem:regdef}.
\end{Lemma}

\begin{Proof}
  If the entries of $\phi$ fail to span $\langle x_1, \ldots, x_5
  \rangle$ then $\Phi$ is clearly not minimal and $v(\Delta(\Phi)) \ge
  12$. Therefore an $\OK$-equivalence brings us to the situation
  considered in Lemma~\ref{lem:regdef}.  Let $d$ be the dimension of
  the tangent space to $C_\phi$ at $P$.  If $d=1$ we are done.  If $d
  \ge 3$ we may assume $\phi_{34} \in \langle x_5 \rangle$ and
  $\phi_{35} = \phi_{45} = 0$. Then
  \[ [\Diag(\pi^{1/2},\pi^{1/2},1,1,1),\pi^{-1/2}
  \Diag(\pi^{-1/2},1,1,1, \pi^{1/2}) ] \Phi \] has coefficients in
  $\OK[\pi^{1/2}]$. So in this case $v(\Delta(\Phi)) \ge 6$.

  Now suppose $d=2$. We may assume $\phi_{34} = x_4$, $\phi_{35} =
  x_5$ and $\phi_{45} = 0$.  To complete the proof we show that if
  $\Phi_{45}$ has coefficient of $x_1$ divisible by $\pi^2$ then
  $v(\Delta(\Phi)) \ge 2$.  Checking this directly, using the formulae
  for the invariants in Section~\ref{sec:inv}, is unfortunately not
  practical.  Instead we argue as follows.  By making substitutions of
  the form $x_4 \leftarrow x_4 + \la x_1$ and $x_5 \leftarrow x_5 +
  \mu x_1$ for suitable $\la,\mu \in \pi \OK$ we may arrange that
  $\Phi_{34}$ and $\Phi_{35}$ also have their coefficients of $x_1$
  divisible by $\pi^2$. Then substituting for $x_1$ we have
  \begin{equation*}
    \Phi = \begin{pmatrix}
      0 & x_1 & \alpha_1 & \alpha_2 & \alpha_3 \\
      & 0 & \beta_1 & \beta_2 & \beta_3 \\
      &   & 0 & \ell_3 & -\ell_2  \\
      &  - &   &  0 & \ell_1 \\
      & & & & 0
    \end{pmatrix}
  \end{equation*}
  where $\ell_1 \equiv 0 \pmod{\pi}$, the coefficient of $x_1$ in each
  of the $\alpha_i$ and $\beta_i$ vanishes mod $\pi$, and the
  coefficient of $x_1$ in each of the $\ell_i$ vanishes mod $\pi^2$.
  By subtracting suitable multiples of the first two rows/columns from
  the last three rows/columns we may further assume that the
  coefficient of $x_1$ in each of the $\alpha_i$ and $\beta_i$
  vanishes mod $\pi^2$.  Since it only matters what the coefficients
  are mod $\pi^2$, we may now assume that none of the $\alpha_i$,
  $\beta_i$ and $\ell_i$ involve $x_1$. By \cite[Lemma 2.4]{minred5},
  $\Phi$ has the same discriminant as the quadric intersection
  \begin{align*}
    \ell_1 \alpha_1 + \ell_2 \alpha_2 + \ell_3 \alpha_3 &= 0 \\
    \ell_1 \beta_1 + \ell_2 \beta_2 + \ell_3 \beta_3 &= 0.
  \end{align*}
  Since $\ell_1 \equiv 0 \pmod{\pi}$, the reduction of this quadric
  intersection mod $\pi$ contains a line. It can then be checked (for
  example by a brute force calculation)
that the discriminant vanishes mod $\pi^2$. This completes the proof.
\end{Proof}


\section{Geometric minimality and an application}
\label{sec:geom-min}

In this section we define the notion of {\em geometric minimality} and
explain the role it has to play in the proof of Theorem~\ref{thm1}.
We assume from now on that the residue field $k$ is algebraically
closed.  Following \cite[Definition 8.3.1]{Liubook} we have
\begin{Definition}
  A {\em fibred surface} $\CC/S$ is an integral projective flat
  $S$-scheme of dimension $2$.
\end{Definition}

\begin{Lemma}
\label{lemma:normal}
Let $C \subset \PP^{n-1}_K$ be a smooth projective curve and $\CC$ its
closure in $\PP^{n-1}_S$. Then $\CC$ is a fibred surface. Moreover
$\CC$ is normal if and only if
\begin{enumerate}
\item $\CC$ is Cohen-Macaulay, and
\item there are only finitely many non-regular points on the special fibre.
\end{enumerate}
\end{Lemma}
\begin{Proof} The coordinate ring of $\CC$ is a subring of that of
  $C$. Since $C$ is integral it follows that $\CC$ is integral.  Then
  $\CC \to S$ is flat and $\dim \CC = 2$ by \cite[Corollaries 4.3.10
  and 4.3.14]{Liubook}.  By definition $\CC$ is projective.  Since
  $\dim \CC=2$ and the generic fibre is smooth, (i) and (ii) are
  equivalent to the conditions $(S_2)$ and $(R_1)$ in Serre's
  criterion \cite[Theorem 8.2.23]{Liubook}.
\end{Proof}

Let $\CC/S$ be a fibred surface. Lipman~\cite{Artin} showed that if
$K$ is complete then $\CC$ admits a {\em desingularisation}
(i.e. resolution of singularities).  If $\CC$ has smooth generic fibre
then the hypothesis that $K$ is complete may be removed, as described
in \cite[Corollary 8.3.51]{Liubook}. If in addition $\CC$ is normal
then by \cite[Proposition 9.3.32]{Liubook} it admits a {\em minimal
  desingularisation}.

\begin{Definition}
\label{def:geom-min}
Let $C \subset \PP^{n-1}_K$ be a genus one normal curve of degree $n$,
with Jacobian $E$.  Let $\CC$ be the closure of $C$ in $\PP^{n-1}_S$.
We say that $C$ is {\em geometrically minimal} if $\CC$ is normal, and
the minimal desingularisation of $\CC$ is isomorphic (as an
$S$-scheme) to the minimal proper regular model of $E$.
\end{Definition}

This definition is {\em not} invariant under changes of co-ordinates
defined over $K$. We remark that if $C$ is geometrically minimal then
$C(K) \not= \emptyset$, and $\CC$ is obtained from the minimal proper
regular model of $E$ by contracting some of the irreducible components
of the special fibre.

Before explaining how geometric minimality is used in the proof of
Theorem~\ref{thm1}, we quote the following lemma.

\begin{Lemma}
\label{lem:h0}
Let $\CC$ be a projective $S$-scheme, and $\LL$ an invertible sheaf
on~$\CC$.
\begin{enumerate}
\item The natural map $H^0(\CC,\LL) \otimes_{\OK} K \to
  H^0(\CC_K,\LL_K)$ is an isomorphism.
\item We have the inequality $\dim_k H^0(\CC_k,\LL_k) \ge \dim_K
  H^0(\CC_K,\LL_K)$.
\item If equality holds in (ii) then $H^0(\CC,\LL)$ is a free
  $\OK$-module and the natural map $H^0(\CC,\LL) \otimes_{\OK} k \to
  H^0(\CC_k,\LL_k)$ is an isomorphism.
\end{enumerate}
\end{Lemma}

\begin{Proof}
  Part (i) is \cite[Corollary 5.2.27]{Liubook} with $A = \OK$ and $B =
  K$. The rest is \cite[Lemma 5.2.31 and Theorem 5.3.20]{Liubook}.
\end{Proof}

\begin{Theorem}
\label{thm:unimodular}
Let $C_1 \subset \PP^{n-1}_K$ and $C_2 \subset \PP^{n-1}_K$ be genus
one normal curves of degree~$n$. Suppose that $C_1$ and $C_2$ are
isomorphic via a change of coordinates given by $B \in \GL_n(K)$. If
$C_1$ and $C_2$ are geometrically minimal, and their Jacobian $E$ has
Kodaira symbol $I_0$ or $I_1$, then $B \in K^\times \GL_n(\OK)$.
\end{Theorem}
\begin{Proof}
  Since the Jacobian $E$ has Kodaira symbol $I_0$ or $I_1$ the special
  fibre of $\EE$ (the minimal proper regular model of $E$) is either a
  smooth curve of genus one, or a rational curve with a node.  Let
  $\CC_i$ be the closure of $C_i$ in $\PP^{n-1}_S$.  Then $\CC_i$ is
  obtained from $\EE$ by contracting some of the irreducible
  components of the special fibre. Since $\EE_k$ is irreducible and
  $\CC_{i,k}$ is a curve it follows that $\CC_i \isom \EE$. We now
  write $f_i : \EE \to \PP^{n-1}_S$ for the embedding with image
  $\CC_i$ and let $\LL_i = f_i^* \OO(1)$.

  Since $C_i = \CC_{i,K}$ is a genus one curve of degree $n$ we have
  $\dim_K H^0(E,\LL_{i,K}) = n$. Since $\CC_{i,k}$ is either a genus
  one curve or a rational curve with a node, and it has degree $n$ by
  \cite[Chapter III, Corollary 9.10]{Hartshorne}, we have $\dim_k
  H^0(\EE_k,\LL_{i,k}) = n$. Then Lemma~\ref{lem:h0} shows that
  $H^0(\EE,\LL_i) \isom \OO_K^n$. Our choice of co-ordinates on
  $\PP^{n-1}_S$ corresponds to a choice of bases for $H^0(\EE,\LL_1)$
  and $H^0(\EE,\LL_2)$. By hypothesis $\LL_{1,K} \isom \LL_{2,K}$, and
  the isomorphism $H^0(E,\LL_{1,K}) \isom H^0(E,\LL_{2,K})$ is given,
  relative to our chosen bases, by some scalar multiple of $B$.

  Let $\LL = \LL_1 \otimes \LL_2^{-1}$. By Lemma~\ref{lem:h0}(ii) both
  $\LL_k$ and its dual $\LL_k^{-1}$ have non-zero global sections.
  Since $\EE_k$ is irreducible it follows that $\LL_k$ is trivial.
  Then by Lemma~\ref{lem:h0}(iii) both $\LL$ and $\LL^{-1}$ have
  global sections that are nowhere vanishing on the special fibre.
  Therefore $\LL$ is trivial and so $\LL_1 \isom \LL_2$. Taking global
  sections gives an isomorphism of $\OK$-modules $H^0(\EE,\LL_1) \isom
  H^0(\EE,\LL_2)$. This isomorphism is again given, relative to our
  chosen bases, by a scalar multiple of $B$.  It follows that $B \in
  K^\times \GL_n(\OK)$.
\end{Proof}

\begin{ProofOf}{Theorem~\ref{thm1}}
  We saw in Lemma~\ref{lem:sqrfree}(i) that for $\Phi \in X_5(\OK)$
  with $v(\Delta(\Phi)) \le 1$, the Jacobian of $C_\Phi$ has Kodaira
  symbol $I_0$ or $I_1$.  We are free to replace $K$ by an unramified
  extension. So by~\cite[Theorem 7.1]{minred5} we may assume that
  $C_\Phi(K) \not= \emptyset$ and likewise for $\Phi'$.  In
  Sections~\ref{sec:flat}, \ref{sec:normal} and~\ref{sec:invdiff} we
  show that, since $\Phi$ and $\Phi'$ are minimal, $C_\Phi$ and
  $C_{\Phi'}$ are geometrically minimal.  Theorem~\ref{thm:unimodular}
  then shows that $B \in K^\times \GL_5(\OK)$ and we are done by
  Remark~\ref{rem1}.
\end{ProofOf}


\section{Minimal Pfaffian models are flat}
\label{sec:flat}

Let $\Phi \in X_5(\OK)$ with reduction $\phi \in X_5(k)$.  In this
section we show that if $\Phi$ is minimal then $\CC_\Phi$ is a fibred
surface.

\begin{Lemma} 
\label{lem:flat}
If $\Phi \in X_5(\OK)$ is non-singular then the following are
equivalent.
\begin{enumerate}
\item $\CC_\Phi$ is the closure of $C_\Phi$ in $\PP^4_S$.
\item $\CC_\Phi$ is a fibred surface.
\item $C_\phi$ is a curve.
\end{enumerate}
\end{Lemma}
\begin{Proof}
  (i) $\Rightarrow$ (ii) $\Rightarrow$ (iii). See
  Lemma~\ref{lemma:normal}
  and \cite[Corollary 4.3.14]{Liubook}. \\
  (iii) $\Rightarrow$ (i).  Let $R = k[x_1, \ldots, x_5]$.  With
  notation as in Section~\ref{sec:inv}, there is a complex of graded
  free $R$-modules
  \begin{equation}
\label{mfr}
0 \ra R(-5) \stackrel{\Pf(\phi)^T}{\ra} R(-3)^5 \stackrel{\phi}{\ra}
R(-2)^5 \stackrel{\Pf(\phi)}{\ra} R.
\end{equation}
Since $C_\phi$ is a curve, this complex is exact by the
Buchsbaum-Eisenbud acyclicity criterion \cite[Theorem 20.9]{E}.

Let $\Pf(\Phi) = (p_1, \ldots, p_5)$. Let ${\mathcal I}$ be the ideal
in $\mathcal R = \OK[x_1,\ldots,x_5]$ generated by $p_1, \ldots, p_5$.
We must show that if $f \in \mathcal R$ and $\pi f \in {\mathcal I}$
then $f \in {\mathcal I}$. We write $\pi f = \sum_{i=1}^5 f_i p_i$ for
some $f_1,\ldots,f_5 \in \mathcal R$. Then $\sum_{i=1}^5 f_i p_i
\equiv 0 \pmod{\pi}$. Since~(\ref{mfr}) is exact it follows that $f_i
= \pi g_i + \sum_{j=1}^5 \Phi_{ij} h_j$ for some $g_1, \ldots,
g_5,h_1, \ldots h_5 \in \mathcal R$. Then $\pi f = \sum_{i=1}^5 f_i
p_i = \pi \sum_{i=1}^5 g_i p_i$ and so $f \in {\mathcal I}$ as
required.
\end{Proof}

\begin{Lemma}
\label{lem:i-ii}
If $\Phi \in X_5(\OK)$ is minimal then
\begin{enumerate}
\item The $4 \times 4$ Pfaffians of $\phi$ are linearly independent.
\item The subscheme $C_\phi \subset \PP^4_k$ does not contain a plane.
\item The entries of $\phi$ span $\langle x_1, \ldots,x_5 \rangle$.
\end{enumerate}
\end{Lemma}
\begin{Proof} This is \cite[Lemma 7.8]{minred5}.
\end{Proof}

\begin{Lemma}
\label{lem:curve}
If $\phi \in X_5(k)$ satisfies conditions (i) and (ii) in
Lemma~\ref{lem:i-ii} then $C_\phi$ is a curve.
\end{Lemma}
\begin{Proof}
  By \cite[Lemma 5.8]{g1inv} every irreducible component of $C_\phi$
  has dimension at least $1$. We must show there are no components of
  dimension $2$ or more. Let $\Sing C_\phi$ be the set of points of
  $C_\phi$ with tangent space of dimension at least $2$. This contains
  all components of $C_\phi$ of dimension $2$ or more.  If $\Sing
  C_\phi$ is contained in a line then we are done. So suppose $P_1,
  P_2, P_3 \in \Sing C_\phi$ span a plane $\Pi$. If $C_\phi$ contains
  each of the lines $P_i P_j$ then it must contain $\Pi$, since
  $C_\phi$ is defined by quadrics. But this is impossible by (ii). We
  may therefore suppose $P_1 P_2 \not\subset C_\phi$.

  A change of co-ordinates gives $P_1 = (1:0: \ldots:0)$ and $P_2 =
  (0:1: \ldots:0)$. If we write $\phi = \sum x_i M_i$ then $M_1$ and
  $M_2$ have rank $2$, but their sum has rank $4$. Therefore $\phi$ is
  equivalent to a model with $\phi_{12} = x_1$, $\phi_{34} = x_2$ and
  all other $\phi_{ij}$ (for $i<j$) linear forms in
  $x_3,x_4,x_5$. Since $P_1,P_2 \in \Sing C_\phi$ it follows that
  $\phi_{35}$ and $\phi_{45}$ are linearly dependent, and $\phi_{15}$
  and $\phi_{25}$ are linearly dependent. Therefore the space of
  linear forms spanned by the entries of the last row/column of $\phi$
  has dimension at most~$2$. Replacing $\phi$ by a $k$-equivalent
  model brings us to the case
  \[ \phi = \begin{pmatrix}
    0 & \xi & \alpha_1 & \alpha_2 & \alpha_3 \\
    & 0 & \beta_1 & \beta_2 & \beta_3 \\
    & & 0 & \eta & 0 \\
    & - & & 0 & 0 \\
    & & & & 0
  \end{pmatrix} \] where $\xi,\eta,\alpha_1, \alpha_2, \alpha_3,
  \beta_1, \beta_2, \beta_3$ are linear forms in $x_1, \ldots, x_5$.
  By (i) the linear forms $\alpha_3$ and $\beta_3$ are linearly
  independent, and $\eta \not=0$. Therefore $C_\phi$ is the union of
  \[ \Gamma_2 = \{ \alpha_3 = \beta_3 = \xi \eta - \alpha_1 \beta_2 +
  \alpha_2 \beta_1 = 0 \} \] and
  \[ \Gamma_3 = \left\{ \rank \begin{pmatrix} \alpha_1 & \alpha_2 & \alpha_3 \\
      \beta_1 & \beta_2 & \beta_3 \end{pmatrix} \le 1 \right\} \cap \{
  \eta = 0 \}. \] We may think of $\Gamma_2$ as a degenerate conic,
  and $\Gamma_3$ as a degenerate twisted cubic.  It remains to show
  that these degenerations are still curves. In the case of $\Gamma_2$
  this is clear by (ii). In the case of $\Gamma_3$ we use the
  following lemma.  The conditions of the lemma are satisfied by (i)
  and (ii).
\end{Proof}

\begin{Lemma}
  Let $\psi$ be a $2 \times 3$ matrix of linear forms in $x_1,
  \ldots,x_4$.  Let $\Gamma_3 \subset \PP^3$ be defined by $\rank \psi
  \le 1$. Suppose that
\begin{enumerate}
\item The $2 \times 2$ minors of $\psi$ span a vector space of
  dimension at least $2$.
\item The subscheme $\Gamma_3 \subset \PP^3$ does not contain a plane.
\end{enumerate}
Then $\Gamma_3$ is a curve.
\end{Lemma}
\begin{Proof}
  Since $\Gamma_3$ is defined by quadrics, any irreducible component
  of dimension $2$ would have degree $1$ or $2$. These possibilities
  are ruled out by (ii) and (i).
\end{Proof}

\begin{Theorem}
\label{thm:flat}
If $\Phi \in X_5(\OK)$ is minimal then $\CC_\Phi$ is a fibred surface.
\end{Theorem}
\begin{Proof} This is immediate from Lemmas~\ref{lem:flat},
  \ref{lem:i-ii} and~\ref{lem:curve}.
\end{Proof}

\section{Minimal Pfaffian models are normal}
\label{sec:normal}

We have seen that if $\Phi \in X_5(\OK)$ is minimal then $\CC_\Phi$ is
a fibred surface. In this section we show that $\CC_\Phi$ is normal.
If $v(\Delta(\Phi)) \le 1$ then Lemma~\ref{lem:reg} already shows that
$\CC_\Phi$ is regular, and hence normal. To treat the general case we
check the conditions in Lemma~\ref{lemma:normal}.

\begin{Lemma}
\label{lem:lci}
If $\Phi \in X_5(\OK)$ is minimal then
\begin{enumerate}
\item $\CC_\Phi$ is a local complete intersection,
\item $\CC_\Phi$ is Cohen-Macaulay.
\end{enumerate}
\end{Lemma}
\begin{Proof} 
  (i) Since $\CC_\Phi \subset \PP^4_S$ has codimension $3$ we must
  show it is locally defined by $3$ equations.  Let $\Pf(\Phi) = (p_1,
  \ldots, p_5)$.  Since $\Phi$ is alternating, the relations
  $\sum_{i=1}^5 p_i \Phi_{ij} = 0$ for $j=4,5$, show that the
  intersection $\CC_\Phi \cap \{ \Phi_{45} \not= 0 \}$ is defined by
  $p_1 = p_2 = p_3 = 0$. By Lemma~\ref{lem:i-ii}(iii)
  the affine pieces $\{ \Phi_{ij} \not= 0 \}$ cover $\PP^4_S$. \\
  (ii) This follows from (i) and~\cite[Corollary 8.2.18]{Liubook}.
\end{Proof}

We prepare to check the second condition in Lemma~\ref{lemma:normal}.
Recall that we assume $k$ is algebraically closed.

\begin{Lemma}
\label{thm:mult}
Let $\phi \in X_5(k)$ satisfy the conditions of Lemma~\ref{lem:i-ii}.
Suppose $C_\phi$ has a multiple component $\Gamma$.  Then after
replacing $\phi$ by a $k$-equivalent model, we are in one of the
following cases
\begin{align*}
  (i) & \quad \begin{pmatrix}
    0 & x_1 & x_2 & * & * \\
    & 0 & * & * & 0 \\
    & & 0 & * & 0 \\
    & - & & 0 & x_5 \\
    & & & & 0
  \end{pmatrix} & (ii) & \quad \begin{pmatrix}
    0 & x_1 & x_2 & x_3 & 0 \\
    & 0 & x_3 & x_4 & 0 \\
    & & 0 & 0 & x_4 \\
    & - & & 0 & x_5 \\
    & & & & 0
  \end{pmatrix} 
\end{align*} \begin{align*} (iii) & \quad \begin{pmatrix}
    0 & x_1 & x_2 & * & * \\
    & 0 & * & * & * \\
    & & 0 & * & * \\
    & - & & 0 & 0 \\
    & & & & 0
  \end{pmatrix} & (iv) & \quad \begin{pmatrix}
    0 & 0 & x_1 & x_2 & x_4 \\
    & 0 & x_2 & x_3 & x_5 \\
    & & 0 & x_5 & 0 \\
    & - & & 0 & 0 \\
    & & & & 0
  \end{pmatrix}
\end{align*}
where the entries $*$ are linear forms in $x_3,x_4,x_5$.  Moreover
$\Gamma = \{x_3 = x_4 = x_5 = 0\}$ in cases (i),(ii),(iii), and
$\Gamma = \{x_1 x_3 - x_2^2 = x_4 = x_5 = 0\}$ in case (iv).
\end{Lemma}
\begin{Proof}
  Lemma~\ref{lem:curve} shows that $C_\phi$ is a curve and so the
  complex~(\ref{mfr}) is exact.  From this minimal free resolution we
  compute that $C_\phi$ has Hilbert polynomial
  \[ h(t) = \binom{t+4}{4} - 5 \binom{t+2}{4} + 5 \binom{t+1}{4} -
  \binom{t-1}{4} = 5t. \] In particular $C_\phi \subset \PP^4$ has
  degree $5$. The multiple component $\Gamma$ must therefore be a line
  or a conic.

  \paragraph{Case $\Gamma$ is a line.} We may assume $\Gamma = \{ x_3
  = x_4 = x_5 = 0 \}$. Then $\phi = \sum x_i M_i$ where all linear
  combinations of $M_1$ and $M_2$ have rank at most $2$. By
  hypothesis, $M_1, \ldots, M_5$ are linearly independent. So we are
  either in case (iii), or $\phi$ takes the form
  \[ \begin{pmatrix}
    0 & x_1 & x_2 & * & * \\
    & 0 & * & \alpha & \beta \\
    & & 0 & \gamma & \delta \\
    & - & & 0 & x_5 \\
    & & & & 0
  \end{pmatrix} \] where the entries $\alpha,\beta,\gamma,\delta$ and
  $*$ are linear forms in $x_3,x_4,x_5$. By row and column operations
  (and substitutions for $x_1$ and $x_2$) we may suppose
  $\alpha,\beta,\gamma,\delta$ do not involve $x_5$. We write $\alpha
  = \alpha_3 x_3 + \alpha_4 x_4$ and likewise for
  $\beta,\gamma,\delta$.
  As shown in \cite[Section 4]{minred5}, $\Gamma$ is a multiple
  component if and only if the determinant of
  \[ \begin{pmatrix}
    \gamma_3 s - \alpha_3 t & \gamma_4 s - \alpha_4 t \\
    \delta_3 s - \beta_3 t & \delta_4 s - \beta_4 t
\end{pmatrix} \] vanishes as a polynomial in $s$ and $t$. If the rows
of this matrix are linearly dependent (over $k$) then we may reduce to
case (i). Otherwise the columns are linearly dependent, and we may
reduce to the case $\alpha_3 = \beta_3 = \gamma_3 = \delta_3=0$ yet
$\alpha_4 \delta_4 - \beta_4 \gamma_4 \not= 0$. Since $C_\phi$ does
not contain the plane $\{ x_4 = x_5 = 0 \}$ it follows that
$\phi_{23}$, and at least one of $\phi_{14}$ and $\phi_{15}$, involves
$x_3$. By a substitution for $x_3$ we may assume $\phi_{23} = x_3$. By
row and columns operations (and substitutions for $x_1$ and $x_2$) we
may assume $\phi_{14}$ and $\phi_{15}$ are multiples of
$x_3$. Replacing the 4th and 5th rows/columns by suitable linear
combinations, and likewise for the 2nd and 3rd rows/columns, brings us
to case (ii).

\paragraph{Case $\Gamma$ is a conic.} We may assume $\Gamma = \{ x_1
x_3 -x_2^2 = x_4 = x_5 = 0 \}$. Let $\Pf(\phi) = (p_1,
\ldots,p_5)$. Replacing $\phi$ by an equivalent model we have
$p_i(x_1,x_2,x_3,0,0)=0$ for $i=1,2,3,4$ and $p_5(x_1,x_2,x_3,0,0)=x_1
x_3 - x_2^2$. Since $\Pf(\phi) \phi = 0$ and $C_\phi$ is a curve, we
may suppose the last row/column of $\phi$ has entries $x_4, x_5,
0,0,0$. As shown in \cite[Section~4]{minred5}, $\Gamma$ is a multiple
component if and only if $\phi_{34}(x_1,x_2,x_3,0,0) = 0$. In this
case $\phi$ is equivalent to a model of the form
\[ \begin{pmatrix} 0 & \xi & x_1 + \langle x_4, x_5 \rangle
  & x_2 + \langle x_4, x_5 \rangle & x_4 \\
  & 0 & x_2 + \langle x_4, x_5 \rangle
  & x_3 + \langle x_4, x_5 \rangle & x_5 \\
  & & 0 &  \langle x_4, x_5 \rangle & 0 \\
  & - & & 0 & 0 \\
  & & & & 0
\end{pmatrix} \] where each $\langle x_4, x_5 \rangle$ denotes some
linear combination of $x_4$ and $x_5$. Subtracting multiples of the
last three rows/columns from the first two row/columns we may suppose
$\xi=0$.  Since the $4 \times 4$ Pfaffians of $\phi$ are linearly
independent we cannot have $\phi_{34}=0$. So making substitutions for
$x_4$ and $x_5$ brings us to the case
\[ \begin{pmatrix} 0 & 0 & x_1 + \langle x_4, x_5 \rangle
  & x_2 + \langle x_4, x_5 \rangle & \langle x_4, x_5 \rangle \\
  & 0 & x_2 + \langle x_4, x_5 \rangle
  & x_3 + \langle x_4, x_5 \rangle &  \langle x_4, x_5 \rangle \\
  & & 0 &  x_5  & 0 \\
  & - & & 0 & 0 \\
  & & & & 0
\end{pmatrix}. \]

If $P \in \GL_2(k)$ then
\[ P \begin{pmatrix} x_1 & x_2 \\ x_2 & x_3 \end{pmatrix} P^T
= \begin{pmatrix} x'_1 & x'_2 \\ x'_2 & x'_3 \end{pmatrix} \] where
$x'_1,x'_2,x'_3$ are linear combinations of $x_1,x_2,x_3$.  Acting on
$\phi$ by a matrix of the form $\Diag(P,P,1)$ we may therefore reduce
to the case $\phi_{15}=x_4$ and $\phi_{25} = x_5$.  Subtracting
multiples of the 5th row/column from the 3rd and 4th rows/columns we
may assume $\phi_{14} = \phi_{23}$. Then making substitutions for
$x_1,x_2,x_3$ brings us to case (iv).
\end{Proof}

The following lemma and its proof could also be used to extend the
algorithms for testing local solubility in \cite{FS} to genus one
curves of degree $5$.

\begin{Lemma}
\label{thm:nonreg}
If $\Phi \in X_5(\OK)$ is minimal then each multiple component
$\Gamma$ of the special fibre $C_\phi$ has at most three non-regular
points.
\end{Lemma}
\begin{Proof}
  We split into the four cases in Lemma~\ref{thm:mult}. \\
  (i) We put
  \begin{align*}
    \Phi_{25} & \equiv \pi ( \alpha_1 x_1 + \alpha_2 x_2 + \ldots )
    \pmod{\pi^2} \\
    \Phi_{35} & \equiv \pi ( \beta_1 x_1 + \beta_2 x_2 + \ldots )
    \pmod{\pi^2}
\end{align*}
where $\alpha_1, \alpha_2, \beta_1, \beta_2 \in k$. We find that
$(s:t:0:0:0) \in \Gamma$ is a non-regular point if and only if the
linear form $s \phi_{34} - t \phi_{24}$ vanishes, or
\[ \beta_1 s^2 + ( \alpha_1 - \beta_2) s t - \alpha_2 t^2 = 0. \] If
$\phi_{24} = \phi_{34} = 0$ then $C_\phi$ is not a curve. If the
quadratic form in $s$ and $t$ vanishes identically then, after
subtracting a multiple of the 1st row/column from the 5th row/column,
we may assume $\alpha_1 = \alpha_2 = \beta_1 = \beta_2=0$.  Since
$\phi_{45} = x_5$ we may assume by a substitution for $x_5$ that
$\Phi_{45} = x_5$.  Then the transformation
\[ [\Diag(\pi,1,1,1,\pi^{-1}),\pi^{-1}\Diag(1,1,\pi,\pi,\pi^2)] \]
shows that $\Phi$ is not minimal. \\
(ii) We put
\begin{align*}
  \Phi_{24} & \equiv x_4 +
  \pi ( \alpha_1 x_1 + \alpha_2 x_2 + \ldots ) \pmod{\pi^2} \\
  \Phi_{25} & \equiv \pi ( \beta_1 x_1 + \beta_2 x_2 + \ldots ) \pmod{\pi^2} \\
  \Phi_{34} & \equiv \pi ( \gamma_1 x_1 + \gamma_2 x_2 + \ldots ) \pmod{\pi^2} \\
  \Phi_{35} & \equiv x_4 + \pi ( \delta_1 x_1 + \delta_2 x_2 + \ldots
  ) \pmod{\pi^2}
\end{align*}
We find that $(s:t:0:0:0) \in \Gamma$ is a non-regular point if and
only if
\[ \gamma_1 s^3 + (\alpha_1 - \gamma_2 - \delta_1) s^2 t - (\alpha_2 +
\beta_1 - \delta_2) s t^2 + \beta_2 t^3 = 0. \] Making a substitution
for $x_4$ we may assume $\alpha_2 = \delta_1=0$.  Subtracting a
multiple of the 1st row/column from the 5th row/column we may assume
$\alpha_1 = \delta_2 = 0$. If the cubic form in $s$ and $t$ vanishes
identically then $\beta_1 = \beta_2 = \gamma_1 = \gamma_2=0$.  Since
$\phi_{45} = x_5$ we may assume by a substitution for $x_5$ that
$\Phi_{45} = x_5$.  Then the transformation
\[ [\Diag(\pi,1,1,\pi^{-1},\pi^{-1}),\pi^{-1}\Diag(1,1,\pi,\pi^2,\pi^3)] \]
shows that $\Phi$ is not minimal. \\
(iii) We put
\[\Phi_{45} \equiv \pi ( \alpha_1 x_1 + \alpha_2 x_2 + \ldots ) \pmod{\pi^2}.\]
We find that $(s:t:0:0:0) \in \Gamma$ is a non-regular point if and
only if $s \phi_{34} - t \phi_{24}$ and $s \phi_{35} - t \phi_{25}$
are linearly dependent, or $\alpha_1 s + \alpha_2 t = 0$.  If the
first of these possibilities is true for all $s$ and $t$, then
$C_\phi$ is not a curve. If $\alpha_1 = \alpha_2 = 0$ then the
transformation
\[  [\Diag(1,1,1,\pi^{-1},\pi^{-1}),\Diag(1,1,\pi,\pi,\pi)] \]
shows that $\Phi$ is not minimal. \\
(iv) We put
\begin{align*}
\Phi_{35} & \equiv \pi ( \alpha_1 x_1 + \alpha_2 x_2
   + \alpha_3 x_3 + \ldots ) \pmod{\pi^2} \\
\Phi_{45} & \equiv \pi ( \beta_1 x_1 + \beta_2 x_2
   + \beta_3 x_3 + \ldots ) \pmod{\pi^2}
\end{align*}
We find that $(s^2:st:t^2:0:0) \in \Gamma$ is a non-regular point
if and only if
\[  \beta_1 s^3 + ( \beta_2 - \alpha_1) s^2 t + (\beta_3 - \alpha_2) st^2
  - \alpha_3 t^3 = 0. \]
Subtracting a multiple of the first two rows/columns from the last
row/column we may assume $\alpha_1 = \alpha_2 = 0$. If the cubic
form in $s$ and $t$ vanishes identically then $\alpha_i = \beta_i =0$
for $i=1,2,3$. Then the transformation
\[  [\Diag(\pi,\pi,1,1,\pi^{-1}),\pi^{-1}\Diag(1,1,1,\pi,\pi)] \]
shows that $\Phi$ is not minimal.
\end{Proof}

\begin{Theorem}
\label{thm:normal}
If $\Phi\in X_5(\OK)$ is minimal then $\CC_{\Phi}$ is a normal fibred
surface.
\end{Theorem}
\begin{Proof}
  In Section~\ref{sec:flat} we showed that $\CC_{\Phi} \subset
  \PP^4_S$ is the closure of $C_\Phi$ and hence a fibred surface.  The
  conditions for normality in Lemma~\ref{lemma:normal} were checked in
  Lemmas~\ref{lem:lci} and~\ref{thm:nonreg}.
\end{Proof}

\section{Minimal Pfaffian models are geometrically minimal}
\label{sec:invdiff}

In this section we show that if $\Phi \in X_5(\OK)$ is minimal and
$C_\Phi(K) \not= \emptyset$ then $C_{\Phi} \subset \PP_K^4$ is
geometrically minimal.  This extends \cite[Theorem~4.1]{minimal} from
genus one curves of degrees $2$, $3$ and $4$, to degree $5$, and could
also be used to prove results analogous to those in \cite{counting}.

\begin{Definition} Let $E/K$ be an elliptic curve with minimal
  Weierstrass equation $y^2 + a_1 xy + a_3 y = x^3 + a_2 x^2 + a_4 x +
  a_6$.  Then
  \[ \omega_E = \frac{dx}{2y + a_1 x + a_3} \in
  H^0(K,\Omega^1_{E/K}) \] is called a {\em N\'eron differential} on
  $E$. It is uniquely determined up to multiplication by elements of
  $\OO_K^\times$.
\end{Definition}

Let $\CC$ be a fibred surface over $S = \Spec \OK$. If $\CC/S$ is a
local complete intersection then we can define the canonical sheaf
$\omega_{\CC/S}$ as in \cite[Definition 6.4.7]{Liubook}. This is an
invertible sheaf on $\CC$. If $\CC$ has generic fibre $E$ then
$H^0(\CC,\omega_{\CC/S})$ is a sub-$\OK$-module of the $1$-dimensional
$K$-vector space $H^0(E,\Omega^1_{E/K})$.

The following theorem and its proof is closely related to
\cite[Theorem 9.4.35]{Liubook}. See also \cite{Conrad}.

\begin{Theorem}
\label{thm:nerondiff}
Let $E/K$ be an elliptic curve, with minimal proper regular model
$\EE/S$. Let $\CC/S$ be a normal fibred surface with generic fibre
isomorphic to $E$, and minimal desingularisation $\Cdes$.  Suppose
$\CC$ is a local complete intersection and $\omega_{\CC/S} = \omega
\mathcal{O}_{\CC}$ for some $\omega \in H^0(E,\Omega^1_{E/K})$.  The
following are equivalent.
\begin{enumerate}
\item We have $\omega_{\CC/S} = \omega_E \mathcal{O}_{\CC}$ where
  $\omega_E$ is a N\'eron differential on $E$.
\item The morphism $\Cdes \to \EE$ (which exists by definition of
  $\EE$) is an isomorphism.
\end{enumerate}
\end{Theorem}

\begin{Proof} 
  (i) $\Rightarrow$ (ii).  Let $f : \Cdes \to \EE$ be the morphism in
  (ii) and $g : \Cdes \to \CC$ the minimal desingularisation. We are
  assuming that $\omega_{\CC/S} = \omega_E \mathcal{O}_{\CC}$, whereas
  \cite[Theorem 9.4.35]{Liubook} gives that $\omega_{\EE/S} = \omega_E
  \mathcal{O}_{\EE}$. Therefore
\begin{equation}
  \label{eqncan1}
f^* \omega_{\EE/S} \isom \omega_E \mathcal{O}_{\Cdes} \isom g^* \omega_{\CC/S}.
\end{equation}

Let $\Gamma$ be an exceptional divisor (or $(-1)$-curve) on $\Cdes$.
Since the desingularisation $g : \Cdes \to \CC$ is minimal, it does
not contract $\Gamma$. Therefore
\begin{equation*}
\omega_{\Cdes/S}|_\Gamma \isom   g^* \omega_{\CC/S} |_\Gamma 
\end{equation*}
By \cite[Corollary 9.3.27]{Liubook} we know that $\omega_{\EE/S}$ is 
globally free. Therefore each of the sheaves in~(\ref{eqncan1}) is globally
free. Writing $K_{\Cdes/S}$ for a canonical divisor on $\Cdes/S$ we have
\[  K_{\Cdes/S} \cdot \Gamma = \deg( \omega_{\Cdes/S}|_\Gamma ) = 
\deg(g^* \omega_{\CC/S} |_\Gamma ) = 0. \]
On the other hand \cite[Proposition 9.3.10]{Liubook} shows that
$K_{\Cdes/S} \cdot \Gamma < 0$. This is a contradiction. 
We deduce that $\Cdes$ has no exceptional divisors. It follows by the 
factorisation theorem \cite[Theorem 9.2.2]{Liubook} that 
$f : \Cdes \to \EE$ is an isomorphism. 

\noindent
(ii) $\Rightarrow$ (i).
Let $F$ be the exceptional locus of the minimal desingularisation
$g : \Cdes \to \CC$. Then 
\begin{equation}
\label{eqn3}
H^0(\Cdes,\omega_{\Cdes/S}) \subset H^0(\Cdes \setminus F,\omega_{\Cdes/S})
 = H^0(\CC \setminus g(F),\omega_{\CC/S}) = H^0(\CC,\omega_{\CC/S}) 
\end{equation}
where the last equality uses that $\CC$ is normal: see 
\cite[Lemma 9.2.17]{Liubook}. We are assuming that $\Cdes \isom \EE$ 
and $\omega_{\CC/S} = \omega \OO_{\CC}$. Therefore $H^0(\Cdes,\omega_{\Cdes/S})
= \omega_E \OK$ and $H^0(\CC,\omega_{\CC/S}) = \omega \OK$. The 
inclusion~(\ref{eqn3}) shows that $\omega_E = h \omega$ for some
$h \in \OK$. 

Since the sheaves $\omega_{\Cdes/S}$ and $g^* \omega_{\CC/S}$
are identical on $\Cdes \setminus F$, the divisor $\divv(h)$ 
on $\Cdes \isom \EE$ is a sum of irreducible components of $F$. 
On the other hand, $\divv(h)$ is a multiple of the special fibre. Since
not all of the irreducible components of the special fibre are contracted,
it follows that $h \in \OO_K^\times$ as required.
\end{Proof}

\begin{Remark}
  Following the proof of \cite[Theorem 9.4.35(a)]{Liubook}, we have
  the following alternative proof of ``(ii) $\Rightarrow$ (i)''.  Let
  $\Gamma_1, \ldots, \Gamma_r$ be the irreducible components of the
  special fibre that are contracted by $g : \Cdes \to \CC$.  By
  \cite[Theorem 9.1.27]{Liubook} the intersection matrix $(\Gamma_i
  \cdot \Gamma_j)$ 
  is negative definite. Since $\Cdes \isom \EE$ is minimal, an
  argument using Castelnuovo's criterion and the adjunction formula
  (see \cite[Example 9.4.19]{Liubook} or \cite[Chapter IV, Theorem
  8.1(b)]{Silverman}) shows that $K_{\Cdes/S} \cdot \Gamma_i = 0$ for
  all $i$. 
  Therefore the contraction morphism $g : \Cdes \to \CC$ satisfies the
  hypotheses of \cite[Corollary 9.4.18]{Liubook}. As a consequence
  $g_* \omega_{\Cdes/S} = \omega_{\CC/S}$ and $g^* \omega_{\CC/S} =
  \omega_{\Cdes/S}$. Therefore $H^0(\CC,\omega_{\CC/S}) =
  H^0(\Cdes,\omega_{\Cdes/S}) = H^0(\EE,\omega_{\EE/S}) = \omega_E
  \OK$.
\end{Remark}

\begin{Theorem}
\label{thm:can}
Let $\Phi \in X_5(\OK)$ be non-singular with reduction $\phi \in
X_5(k)$.  Suppose $\CC = \CC_\Phi$ is a fibred surface, and the
entries of $\phi$ span $\langle x_1,\ldots,x_5 \rangle$. Then $\CC$ is
a local complete intersection with $\omega_{\CC/S} = \omega_\Phi
\OO_{\CC}$ where $\omega_\Phi$ is defined by~(\ref{invdiff}).
\end{Theorem}

\begin{Proof}
  Exactly as in the proof of Lemma~\ref{lem:lci}, the affine piece
  $\CC \cap \{ \Phi_{45} \not = 0 \}$ is defined by $p_1 = p_2 = p_3 =
  0$. The restriction of the canonical sheaf to this affine piece is
  as claimed by \cite[Corollary 6.4.14]{Liubook} and the next lemma.
  Since the definition~(\ref{invdiff}) of $\omega_\Phi$ is invariant
  under all even permutations of the subscripts, and the affine pieces
  $\{ \Phi_{ij} \not= 0 \}$ cover $\PP^4_S$, the theorem follows.
\end{Proof}

\begin{Lemma}
\label{lem:determinant}
Let $R$ be any ring. Let $\phi \in X_5(R)$ with $\Pf(\phi) = (p_1,
\ldots, p_5)$.  Let $I$ be the ideal in $R[x_1, \ldots,x_5]$ generated
by $p_1, \ldots, p_5$. Then
\begin{equation}
\label{myidentity}
 \frac{\partial(p_1,p_2,p_3)}{\partial(x_1,x_2,x_3)}
  \equiv \phi_{45} \sum_{i,j=1}^5 \frac{\partial p_i}{\partial x_1}
\frac{\partial \phi_{ij}}{\partial x_2}\frac{\partial p_j}{\partial x_3}
\pmod{I}.
\end{equation}
\end{Lemma}

\begin{Proof}
  We have $\sum_{i=1}^5 p_i \phi_{ij} = 0$ for all $1 \le j \le
  5$. Differentiating with respect to $x_k$, and working mod $I$, this
  gives
\begin{equation}
\label{myreln}
\sum_{i=1}^5 \frac{\partial p_i}{\partial x_k} \phi_{ij} = 0.
\end{equation}
Using first that $\phi$ is alternating, and then~(\ref{myreln}), we
compute
\begin{align*}
\phi_{45} \sum_{i=1}^3 \sum_{j=4}^5 \frac{\partial p_i}{\partial x_1}
\frac{\partial \phi_{ij}}{\partial x_2}\frac{\partial p_j}{\partial x_3}
& = \sum_{i=1}^3 \sum_{j=4}^5  \frac{\partial p_i}{\partial x_1}
\left( \frac{\partial \phi_{i4}}{\partial x_2} \phi_{j5}
  -  \frac{\partial \phi_{i5}}{\partial x_2} \phi_{j4} \right)
\frac{\partial p_j}{\partial x_3} \\
& = - \sum_{i,j=1}^3  \frac{\partial p_i}{\partial x_1}
\left( \frac{\partial \phi_{i4}}{\partial x_2} \phi_{j5}
  -  \frac{\partial \phi_{i5}}{\partial x_2} \phi_{j4} \right)
\frac{\partial p_j}{\partial x_3}.
\end{align*}
Subtracting the same identity with $\frac{\partial}{\partial x_1}$
and $\frac{\partial}{\partial x_3}$ switched gives
\begin{equation}
\label{eqn1}
\phi_{45} \sum_{i=1}^3 \sum_{j=4}^5 \frac{\partial \phi_{ij}}{\partial x_2}
 \frac{\partial(p_i,p_j)}{\partial(x_1,x_3)}
= \sum_{i<j} \frac{\partial}{\partial x_2} ( - \phi_{i4} \phi_{j5}
 + \phi_{i5} \phi_{j4} ) \frac{\partial(p_i,p_j)}{\partial(x_1,x_3)}
\end{equation}
where we write $\sum_{i<j}$ for the sum over all $1 \le i < j \le 3$.
Again using~(\ref{myreln}),
\[ \sum_{i,j=4}^5 \frac{\partial p_i}{\partial x_1}
\phi_{ij} \frac{\partial p_j}{\partial x_3}   =
- \sum_{i=4}^5 \sum_{j=1}^3 \frac{\partial p_i}{\partial x_1}
\phi_{ij} \frac{\partial p_j}{\partial x_3}  =
\sum_{i,j=1}^3 \frac{\partial p_i}{\partial x_1}
\phi_{ij} \frac{\partial p_j}{\partial x_3}. \]
Therefore
\begin{equation}
\label{eqn2}
\phi_{45} \frac{\partial(p_4,p_5)}{\partial(x_1,x_3)}
= \sum_{i<j} \phi_{ij} \frac{\partial(p_i,p_j)}{\partial(x_1,x_3)}.
\end{equation}
We break up the sum on the right of~(\ref{myidentity}) as
\[ \sum_{i<j} \frac{\partial \phi_{ij}}{\partial x_2}
\frac{\partial(p_i,p_j)}{\partial(x_1,x_3)}
+ \sum_{i=1}^3 \sum_{j=4}^5 \frac{\partial \phi_{ij}}{\partial x_2}
\frac{\partial(p_i,p_j)}{\partial(x_1,x_3)} +
\frac{\partial \phi_{45}}{\partial x_2}
\frac{\partial(p_4,p_5)}{\partial(x_1,x_3)}.   \]
Then by~(\ref{eqn1}) and~(\ref{eqn2}), the right hand side
of~(\ref{myidentity}) is
\[ \sum_{i<j} \frac{\partial}{\partial x_2} (\phi_{ij} \phi_{45}
 - \phi_{i4} \phi_{j5} + \phi_{i5} \phi_{j4} )
\frac{\partial(p_i,p_j)}{\partial(x_1,x_3)}. \]
Since for $i,j,k$ an even permutation of $1,2,3$ we have
$-p_k = \phi_{ij} \phi_{45}
 - \phi_{i4} \phi_{j5} + \phi_{i5} \phi_{j4}$
the result follows.
\end{Proof}

\begin{Remark} 
\label{char2}
We keep the notation of the lemma.  Differentiating the relation
$\sum_{j=1}^5 \phi_{ij} p_j = 0$ with respect to $x_4$ and $x_5$ we
have
\[ \sum_{j=1}^5 \frac{\partial \phi_{ij}}{\partial x_4}
   \frac{\partial p_j}{\partial x_5}  
+ \sum_{j=1}^5 \frac{\partial \phi_{ij}}{\partial x_5}
   \frac{\partial p_j}{\partial x_4}  
+ \sum_{j=1}^5 \phi_{ij}
   \frac{\partial^2 p_j}{\partial x_4 \partial x_5} = 0. \]
We multiply by $\frac{\partial p_i}{\partial x_4}$ and sum over $i$.
By~(\ref{myreln}) and the fact $\phi$ is alternating, 
the second two terms vanish mod $I$. Therefore 
\[ \sum_{i,j=1}^5 \frac{\partial p_i}{\partial x_4}
  \frac{\partial \phi_{ij}}{\partial x_4}
   \frac{\partial p_j}{\partial x_5} \equiv 0 \pmod{I}. \]
This shows that the restriction $\Char K \not= 2$ in 
\cite[Section 5.4]{g1inv} is not needed.
\end{Remark}

\begin{Lemma}
\label{lem:Level0Neron}
Let $\Phi \in X_5(K)$ be non-singular with $C_\Phi(K) \not=
\emptyset$.  Then $\Phi$ has level~$0$ if and only if $\omega_\Phi$ is
a N\'eron differential on $C_\Phi \isom E$.
\end{Lemma}

\begin{Proof}
Let $E/K$ have minimal Weierstrass equation 
\[ y^2 + a_1 xy + a_3 y = x^3 + a_2 x^2 + a_4 x + a_6. \] The complete
linear system $|4.0_E|$ defines a morphism $\alpha : E \to \PP^3$. It
is given by $(x,y) \mapsto (1:x:y:x^2)$. The image is $C_4 = \{ Q_1 =
Q_2 = 0 \} \subset \PP^3$ where
\begin{align*}
Q_1 & = x_1 x_4 - x_2^2, \\
Q_2 & = x_3^2  + a_1 x_2 x_3 + a_3 x_1 x_3 - x_2 x_4 - a_2 x_2^2 
- a_4 x_1 x_2 - a_6 x_1^2,
\end{align*}
and an invariant differential $\omega_4$ on $C_4$ is given by
\[ \omega_4 = \frac{x_1^2 d(x_2/x_1)}{ \frac{\partial Q_1}{\partial
    x_4} \frac{\partial Q_2}{\partial x_3} - \frac{\partial
    Q_1}{\partial x_3} \frac{\partial Q_2}{\partial x_4}}. \] We claim
that (i) $\Delta(Q_1,Q_2) = \Delta_E$ and (ii) $\omega_4$ is a N\'eron
differential on $C_4 \isom E$. Indeed the invariants were scaled in
\cite{g1inv} so that (i) is true, whereas for (ii) it is easy to see 
that $\alpha^* \omega_4 = dx/(2y + a_1 x + a_3)$.

Since $C_\Phi(K) \not= \emptyset$ we may identify $C_\Phi \isom E$.
The hyperplane section is linearly equivalent to $4.0_E + P$ for some
$P \in E(K)$.  Let $\Psi \in X_5(K)$ be the Pfaffian model constructed
from the quadric intersection $(Q_1,Q_2)$ by ``unprojection centred at
$P$'' as described in \cite[Lemma 2.3]{minred5}.  By \cite[Lemma
2.4]{minred5} and its proof, we have (i) $\Delta(\Psi) = \Delta_E$ and
(ii) $\omega_\Psi$ is a N\'eron differential on $C_\Psi \isom E$.

The curves $C_\Phi$ and $C_\Psi$ differ by a change of co-ordinates
defined over $K$. So by \cite[Theorem 4.1(ii)]{explicit5}, 
the Pfaffian models $\Phi$ and $\Psi$ are $K$-equivalent, say $\Phi = 
g \Psi$ for some $g \in \GL_5(K) \times \GL_5(K)$. Since $\Delta$
is an invariant of weight $12$ we have $\Delta(\Phi)
 = (\det g)^{12} \Delta(\Psi)$. 
Let $\gamma : C_\Phi \to C_\Psi$ be the isomorphism described by $g$. 
By \cite[Proposition 5.19]{g1inv} we have $\gamma^* \omega_\Psi = 
(\det g) \omega_\Phi$. Therefore both the conditions in the 
statement of the lemma are equivalent to $v(\det g)=0$.
\end{Proof}

\begin{Remark}
If $\Char K \not=2,3$ then \cite[Proposition 5.23]{g1inv}
shows that $(C_\Phi,\omega_\Phi)$ and $(E, \omega)$ are isomorphic
over $\Kbar$, where $E$ is the elliptic curve~(\ref{eqn:jac})
and $\omega = dx/(2y + a_1(\Phi) x + a_3(\Phi))$.  
This gives an alternative proof of Lemma~\ref{lem:Level0Neron}.
The isomorphism $C_\Phi \isom E$ might not be defined over $K$,
but differs from an isomorphism that is defined over $K$ by
an automorphism of the curve $E$. The latter might rescale 
$\omega$ by a root of unity, but won't change whether it is a
N\'eron differential.
\end{Remark}



\begin{Theorem}
\label{thm:minequiv}
Let $\Phi \in X_5(\OK)$ be non-singular with $C_\Phi(K) \not=
\emptyset$.  Suppose $\CC_\Phi$ is a fibred surface, and the entries
of $\phi$ span $\langle x_1,\ldots,x_5 \rangle$.  Then $\Phi$ is
minimal if and only if $C_\Phi$ is geometrically minimal.
\end{Theorem}
\begin{Proof} 
  Lemma~\ref{lem:flat} shows that $\CC = \CC_\Phi$ is the closure of
  $C_\Phi$ in $\PP^4_S$. By either Definition~\ref{def:geom-min} or
  Theorem~\ref{thm:normal} we may suppose $\CC$ is normal.  Let $E$ be
  the Jacobian of $C_\Phi$.  Since $C_\Phi(K) \not= \emptyset$ we have
  $C_\Phi \isom E$.  Theorem~\ref{thm2} and
  Lemma~\ref{lem:Level0Neron} show that $\Phi$ is minimal if and only
  if $\omega_\Phi$ is a N\'eron differential on $C_\Phi \isom E$.  The
  theorem now follows from Theorems~\ref{thm:nerondiff}
  and~\ref{thm:can}.
\end{Proof}

By Lemma~\ref{lem:i-ii}, Theorem~\ref{thm:flat} and
Theorem~\ref{thm:minequiv} we have
\begin{Corollary}
  If $\Phi \in X_5(\OK)$ is minimal and $C_\Phi(K) \not= \emptyset$
  then $C_{\Phi}$ is geometrically minimal.
\end{Corollary}


\section{An alternative proof of Theorem~\ref{thm1}}
\label{sec:alt}

We give a short alternative proof of Theorem~\ref{thm1}, that avoids
using schemes, except for the definition of a regular point.  It would
however be rather hard to motivate this proof without the work in
earlier sections.

By putting the matrices $A, B \in \GL_5(K)$ in Smith normal form
(and making use of Remark~\ref{rem1}),
Theorem~\ref{thm1} is equivalent to the following.

\begin{Theorem}
\label{thm1alt}
Let $\Phi, \Phi' \in X_5(\OK)$ with $v(\Delta(\Phi))
\le 1$ and $v(\Delta(\Phi')) \le 1$. If
  \begin{equation*}
  \Phi' =  [\Diag(\pi^{-r_1}, \ldots, \pi^{-r_5}),\Diag(\pi^{s_1}, \ldots,
    \pi^{s_5}) ] \Phi
\end{equation*}
for some $r_1, \ldots, r_5, s_1, \ldots, s_5 \in \Z$ then
$s_1 = s_2 = \ldots = s_5$.
\end{Theorem}

For the proof we may assume the residue field $k$ is algebraically
closed. As before we write $\phi \in X_5(k)$ for the reduction of
$\Phi$ mod $\pi$.  For the purposes of this section, a $k$-point $P$
on $C_\phi$ is {\em regular} if it satisfies the conditions in
Lemma~\ref{lem:regdef}, and otherwise {\em non-regular}.  Since $\dim
\CC_\Phi = 2$ this agrees with the standard terminology, but we don't
need to know this.

\begin{Lemma} If $v(\Delta(\Phi)) \le 1$ then $C_\phi$ contains no
  lines or conics.
\end{Lemma}

\begin{Proof}
If $C_\phi$ contains a line or conic then, arguing as in the proof of
Lemma~\ref{thm:mult}, we may assume
\begin{equation*}
\phi = \begin{pmatrix}
  0 & x_1 & x_2 & * & * \\
  & 0 & * & * & * \\
  &   & 0 & * & *  \\
  &  - &   &  0 & * \\
  & & & & 0
\end{pmatrix}
\quad \text{ or } \quad
\begin{pmatrix}
  0 & * & * & * & * \\
  & 0 & * & * & * \\
  &   & 0 & * & 0  \\
  &  - &   &  0 & 0 \\
  & & & & 0
\end{pmatrix}
\end{equation*}
where the entries $*$ on the left are linear forms in $x_3,x_4,x_5$,
and on the right are linear forms in $x_1, \ldots, x_5$.
In the first case we apply the transformation
  \[ [\Diag(\pi,1,1,1,1),\pi^{-1} \Diag(1,1,\pi,\pi,\pi) ]. \]
Then $\phi_{14} = \phi_{15} = 0$ and an $\OK$-equivalence brings
us to the second case. In the second case we may assume
$\phi_{34} \in \langle x_1 \rangle$. Applying the transformation
\[ [\Diag(\pi,\pi,1,1,1),\pi^{-1} \Diag(\pi,1,1,1,1) ] \]
gives a model with a non-regular point at $(1:0:\ldots:0)$.
Since all transformations we have used preserve (the valuation of) the
discriminant, we are done by Lemma~\ref{lem:reg}.
\end{Proof}

\begin{Lemma}
\label{getineq}
Let $\Phi, \Phi' \in X_5(\OK)$ be Pfaffian models satisfying
  \begin{equation*}
  \Phi' =  [\Diag(\pi^{-r_1}, \ldots, \pi^{-r_5}),\Diag(\pi^{s_1}, \ldots,
    \pi^{s_5}) ] \Phi
\end{equation*}
for some $r_1 \le \ldots \le r_5$ and $s_1 \le  \ldots \le s_5$.
\begin{enumerate}
\item
If $C_\phi$ contains no lines then $r_1 + r_4 \le s_2$, $r_2 + r_3
\le s_2$ and $r_2 + r_4 \le s_3$.
\item
If $C_\phi$ contains no lines or conics then $r_1 + r_5 \le s_3$, $r_2 + r_5
\le s_4$, $r_3 + r_4 \le s_4$ and $r_3 + r_5 \le s_5$.
\end{enumerate}
\end{Lemma}
\begin{Proof}
(i) If $r_1 + r_4 > s_2$ then all entries of $\phi$ outside the top
left $3 \times 3$ submatrix are linear forms in $x_3,x_4,x_5$. So
$C_\phi$ contains the line $\{x_3 = x_4 = x_5 =0\}$. If $r_2 + r_3 > s_2$
then all entries of $\phi$ outside the first row/column
are linear forms in $x_3,x_4,x_5$. So
$C_\phi$ contains the line $\{x_3 = x_4 = x_5 =0\}$.
If $r_2 + r_4 > s_3$ then $C_\phi$ contains the line $\{\phi_{23} = x_4 =
x_5  =0 \}$.

(ii) If $r_1 + r_5 > s_3$, $r_2 + r_5 > s_4$ or $r_3 + r_5 > s_5$ then the
entries of the last row/column of $\phi$ are in $\langle x_4, x_5
\rangle$, $\langle \phi_{15}, x_5 \rangle$ or $\langle \phi_{15},
\phi_{25} \rangle$. If $r_3 + r_4 > s_4$ then the bottom
right $3 \times 3$ submatrix of $\phi$ has entries in
$\langle x_5 \rangle$.
So in all these cases
$\phi$ is $k$-equivalent to a model with $\phi_{35} = \phi_{45} = 0$.
Let $p_5$ be the Pfaffian of the top left $4 \times 4$ submatrix.
Then $C_\phi$ contains $\{ \phi_{12} = \phi_{25} = p_5 = 0\}$ which
is either a conic or contains a line.
\end{Proof}

\begin{Lemma}
\label{getap}
Let $\Phi$ and $\Phi'$ be as in Theorem~\ref{thm1alt}, and suppose
$0 = r_1 \le \ldots \le r_5$ and $s_1 \le \ldots \le s_5$. Then
the $r_i$ and $s_i$ are given by
\[\begin{array}{ccccc|ccccc}
r_1 & r_2 & r_3 & r_4 & r_5 & s_1 & s_2 & s_3 & s_4 & s_5 \\ \hline
0 & \alpha & 2 \alpha & 3 \alpha & 4 \alpha & \le2\alpha & 3 \alpha
& 4 \alpha & 5 \alpha & \ge6\alpha
\end{array}\]
for some $\alpha \ge 0$. 
\end{Lemma}

\begin{Proof} The inequalities in Lemma~\ref{getineq} together
with the inequalities obtained when we replace
$(r_1, \ldots, r_5;s_1, \ldots, s_5)$ by
$(-r_5, \ldots, -r_1;-s_5, \ldots, -s_1)$
give
\begin{align*}
s_2 & = r_1 + r_4 = r_2 + r_3 \quad \implies r_2 - r_1 = r_4 - r_3 \\
s_3 & = r_1 + r_5 = r_2 + r_4 \quad \implies r_2 - r_1 = r_5 - r_4 \\
s_4 & = r_2 + r_5 = r_3 + r_4 \quad \implies r_3 - r_2 = r_5 - r_4
\end{align*}
Therefore $r_1, \ldots, r_5$ are in arithmetic progression.
The other statements follow.
\end{Proof}

\begin{Lemma}
\label{getineq2}
Let $\Phi, \Phi' \in X_5(\OK)$ be Pfaffian models satisfying
  \begin{equation*}
  \Phi' =  [\Diag(\pi^{-r_1}, \ldots, \pi^{-r_5}),\Diag(\pi^{s_1}, \ldots,
    \pi^{s_5}) ] \Phi
\end{equation*}
for some $r_1 \le \ldots \le r_5$ and $s_1 \le  \ldots \le s_5$.
\begin{enumerate}
\item If $r_1 + r_4 > s_1$ and $r_4 + r_5 > s_5 > s_1$ then $\CC_\Phi$ has a
non-regular point.
\item If $r_1 + r_3 > s_1$ and $r_3 + r_4 > s_3 > s_1$ then $\CC_\Phi$ has a
non-regular point.
\item If $r_2 + r_5 < s_5$ and $r_1 + r_2 < s_1 < s_5$ then $\CC_{\Phi'}$ has a
non-regular point.
\item If $r_3 + r_5 < s_5$ and $r_2 + r_3 < s_3 < s_5$ then $\CC_{\Phi'}$ has a
non-regular point.
\end{enumerate}
\end{Lemma}
\begin{Proof}
(i) Since $r_1 + r_4 > s_1$ the only entries of $\phi$ involving $x_1$
are in the top left $3 \times 3$ submatrix. So $P = (1:0: \ldots:0)$
is a point on $C_\phi$. Since $r_4 + r_5 > s_5$ we have $\phi_{45} = 0$
and so $P$ is a singular point. Since $r_4 + r_5 > s_1 + 1$ the
coefficient of $x_1$ in $\Phi_{45}$ vanishes mod $\pi^2$. Therefore $P$
is a non-regular point.

(ii) Since $r_1 + r_3 > s_1$ the only entries of $\phi$ involving $x_1$
are in the top left $2 \times 2$ submatrix. So $P = (1:0: \ldots:0)$
is a point on $C_\phi$. Since $r_3 + r_4 > s_3$ we have $\phi_{34}, \phi_{35},
\phi_{45} \in \langle x_4, x_5 \rangle$
and so $P$ is a singular point. Since $r_3 + r_4 > s_1 + 1$ the
coefficient of $x_1$ in each of $\Phi_{34}, \Phi_{35}$ and $\Phi_{45}$
vanishes mod $\pi^2$. Therefore $P$ is a non-regular point.

(iii), (iv) These follow from (i) and (ii) by switching the roles
of $\Phi$ and $\Phi'$.
\end{Proof}

\begin{ProofOf}{Theorem \ref{thm1alt}}
We may assume $r_1 \le \ldots \le r_5$ and $s_1 \le \ldots \le s_5$.
Replacing $r_i$ by $r_i + \lambda$ and $s_i$ by $s_i + 2\lambda$
still gives the same transformation. So we may assume $r_1=0$.
Then the $r_i$ and $s_i$ are as given in Lemma~\ref{getap}.

If $\alpha = 0$ then $r_1 = \ldots = r_5$ and the conclusion
$s_1 = \ldots =s_5$ follows from the fact $\Phi$ and $\Phi'$ are
minimal. We assume for a contradiction that $\alpha \ge 1$.
Since $r_1 + r_4 = 3 \alpha > s_1$ it follows by
Lemmas~\ref{lem:reg} and~\ref{getineq2}(i) that $r_4 + r_5 \le s_5$.
Since $r_2 + r_3 = 3 \alpha < s_3$ it follows by
Lemmas~\ref{lem:reg} and~\ref{getineq2}(iv)
that $r_3 + r_5 \ge s_5$.
Putting these together we have
\[ r_4 + r_5 \le s_5 \le r_3 + r_5. \]
Therefore $r_3 = r_4$ and this contradicts our assumption that
$\alpha \ge 1$.
\end{ProofOf}

\end{document}